\documentclass[11pt]{amsart}

\usepackage{amssymb,amsmath,amsfonts,eurosym,geometry,ulem,graphicx,caption,xcolor,setspace,comment}
 \usepackage{amsaddr}
\usepackage[perc]{overpic}
\usepackage[all]{xy}
\usepackage{hyperref}

\usepackage[
    natbib=true,
    style=alphabetic,
]{biblatex}

\newtheorem{theorem}{Theorem}[section]
\newtheorem{corollary}[theorem]{Corollary}
\newtheorem{proposition}[theorem]{Proposition}
\newtheorem{lemma}[theorem]{Lemma}

\newtheorem{remark}[theorem]{Remark}

\DeclareMathOperator{\R}{\mathbb{R}}

\DeclareMathOperator{\PP}{\mathcal{P}}

\title[Analytical study of The Lorenz system]{Analytical study of The Lorenz system: existence of infinitely many periodic orbits and their topological characterization}
\author{Tali Pinsky}
\address{The Technion, Haifa}
\email{talipi@technion.ac.il}
\thanks{This research was supported by the Israel Science Foundation (grant No. 51/4051).}
\begin{document}
\maketitle

\begin{abstract}
We consider the Lorenz equations, a system of three dimensional ordinary differential equations modeling atmospheric convection. These equations are chaotic and hard to study even numerically, and so a simpler ``geometric model" has been introduced in the seventies. One of the classical problems in dynamical systems is to relate the original equations to the geometric model. This has been achieved numerically by Tucker for the classical parameter values, and remains open for general values. In this paper we establish analytically a relation to the geometric model, for a different set of parameter values that we prove must exist. This is facilitated by finding a novel way to apply topological tools developed for the study of surface dynamics to the more intricate case of three dimensional flows. 
\end{abstract}

\section{Introduction} \label{sec:Introduction}
In 1963 Edward Lorenz, a  mathematician and meteorologist, introduced a simplified model for heat convention~\cite{lorenz1963deterministic}:
\begin{equation}\left\{\begin{array}{rlr}
\dot{x}(t)  & = \sigma(y-x)\\
\dot{y}(t) & =\rho x -y - xz \\
\dot{z}(t) & =xy-\beta z
\end{array}\right.
\end{equation}

This system has been extensively studied over the last decades, being the first system to demonstrate chaos in a deterministic low dimensional setting and serving as a paradigm for chaotic systems \cite{Ghys2013LorenzParadigm}. It arises in many different mathematical models of chaotic systems and has been studied especially at the classical parameter values $\rho=28, \sigma=10, \beta=\frac{8}{3}$ which were those originally studied by Lorenz.

At these parameter values, the three-dimensional motion given by the Lorenz equations (1) converges almost always onto its famous butterfly-shaped strange attractor. Although this is easily seen to be the case by running any ODE solver, there is no analytical proof for it. 

In the seventies, a simpler dynamical system called the geometric Lorenz model, was developed  \cite{1977OriginAndStructure, Guckenheimer1979}. The geometric model shares by design much of the characteristics of the equations, and also displays a butterfly-shaped strange attractor. At the same time, the geometric model is susceptible to analytic study and it can be proven to be chaotic in a specific sense. i.e., one can prove both geometric properties on the shape of the attractor and the orbits within it, and statistical properties about the mixing rate, for example, the attractor has a central limit theorem \cite{HollandMelbourne2007} and exponential decay of correlations \cite{AraujoMelbourne2016}. The knots that appear as periodic orbits in the geometric model, called ``Lorenz knots'', have been studied topologically by Williams \cite{Williams1979} and Birman and Williams \cite{BirmanWilliams:KnottedOrbitsI} (see also  \cite{Dehornoy2011Noeds, BirmanKofmantwist}) and have a number of surprising knot properties. 

In 1998 Steve Smale compiled a list of 18 major unsolved problems in mathematics ~\cite{smaleproblems}, many of which are still unresolved. His 14th problem  asks whether the Lorenz system at the classical parameters given above  can be proven to be equivalent to the geometric model. This has been answered by  Warwick Tucker \cite{Tucker1999}, via a rigorous numerical proof. Tucker's proof implies in particular that the Lorenz equations indeed possess a butterfly strange attractor.

Ghys \cite{Ghys:KnotsDynamics} proved that the Lorenz knots also appear as the set of periodic orbits for a very different flow, the geodesic flow on the modular surface. 

\begin{figure}[ht]
    \centering
    \begin{overpic}[width=9cm]{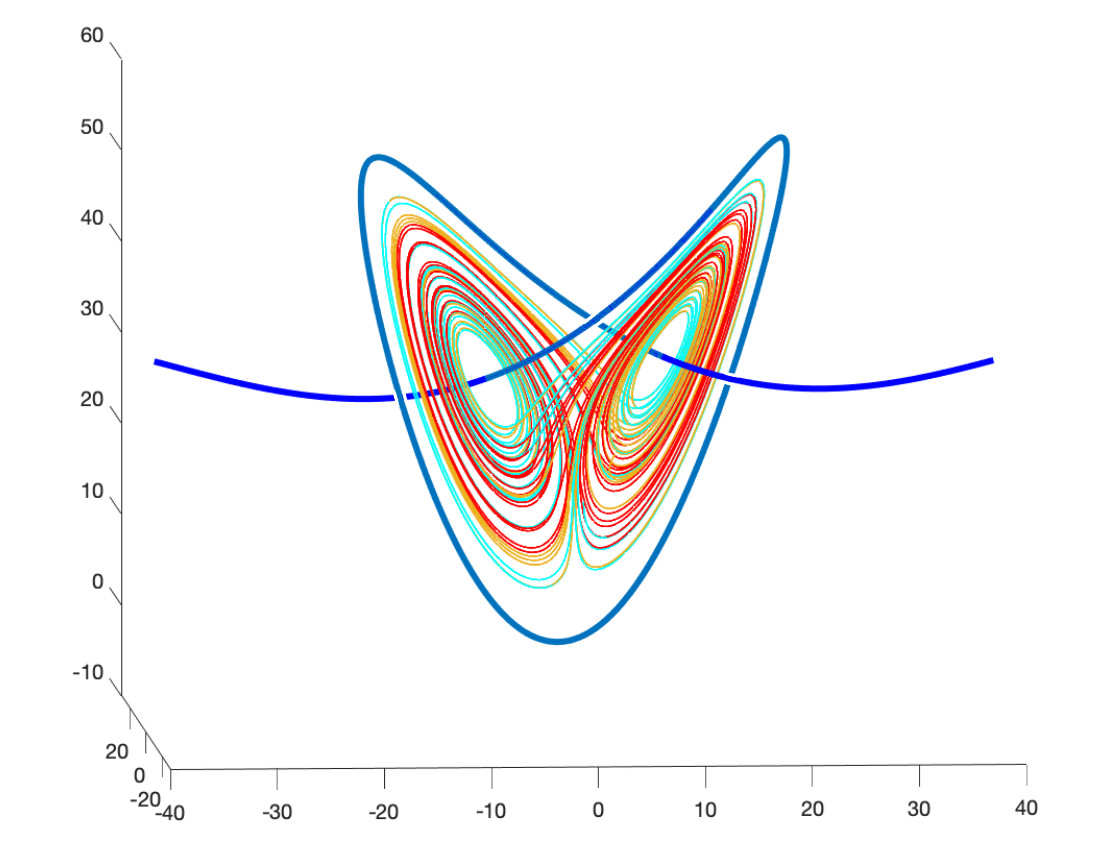}
    \put(1,72){$z$}
    \put(90,0){$x$}
    \end{overpic}
    \caption{The trefoil knot (in bold) at the first T-point, (c.f. \cite{creaser15alpha}).}
    \label{fig:trefoil}
\end{figure}

This is surprising as the modular geodesic flow is a hyperbolic flow arising from number theory. It is chaotic as well, but much better understood and proving its geometric and statistical properties 
is much more approachable. The modular flow is not defined on the entire three dimensional space. One must exclude a one dimensional set that is a knot, called the trefoil knot, which is the thick knot depicted in blue in Figure~\ref{fig:trefoil}. In knot theory, the trefoil is the simplest nontrivial knot.

The Lorenz equations have three singular points where the vector field defining the equations vanishes.
Thus, one may search for ``heteroclinic connections'', i.e. regular orbits that limit onto one of these singular points as $t\to-\infty$ and limit onto another as $t\to\infty$. It is well known numerically that there exist parameters of the Lorenz equations  (called T-points in the literature) for which there is a heteroclinic connection connecting the origin to one of the singularities in the center of the butterfly wings \cite{petrov80homo, Bykov1980, alfsen85system}. Due to the symmetry $(x,y,z)\mapsto (-x,-y,z)$ for the equations, it follows that in these parameters there is also a symmetric heteroclinic connection connecting the origin to the other wing center, and thus all three singular points are connected there.
It has been numerically observed in \cite{pinsky2017topology} that by adding the other side of the stable manifold of each of the the wing centers, this symmetric pair of heteroclinic connections can be continued to an invariant one dimensional knot passing also through infinity. For the first T-point at 
$\rho_0\approx 30.8680, \sigma_0\approx 10.1673$ and $\beta_0=\frac{8}{3}$,
the invariant knot is a trefoil knot, depicted in Figure~\ref{fig:trefoil}.

In an attempt to explain the relation between the equations, the geometric model and the modular geodesic flow, Christian Bonatti and the author have constructed an extension of the geometric model \cite{BonattiPinsky2021}, 
which is  defined on the three-dimensional sphere.
This extended model $X_r$ interpolates between a flow $X_{-1}$ on the three sphere, similar to the classical geometric Lorenz model and having a Lorenz attractor, and a flow $X_1$ that has a hyperbolic non wandering set and an invariant trefoil shaped heteroclinic connection. It is proven in \cite{BonattiPinsky2021} that when considering $X_1$ on the complement of a neighborhood of the trefoil knot, it is topologically equivalent (i.e. there is a homeomorphism taking orbits of one flow to orbits of the other) to the modular geodesic flow, when the latter is compactified as in \cite{Ghys:KnotsDynamics} (turning the one dimensional cusp into a two dimensional torus boundary).

Our main theorem is an approximate answer to Smale's problem. It is an analytical proof that there indeed exists a parameter for which one can define the equations on the trefoil complement, that the equations have a cross section there, and that at this parameter the Lorenz equations can be continuously deformed to the extended geometric model $X_1$ (having a hyperbolic nonwandering set) while preserving the cross section.

\begin{theorem}\label{thm:MainThm}
There exists a parameter value for the Lorenz equations for which a trefoil-shaped heteroclinic connection exists.
Furthermore, there is a one parametric family of flows on $S^3$, deforming continuously the Lorenz equations at the trefoil parameter to the extended geometric model $X_1$. For any flow along the deformation there exists  a cross section, so that the return map corresponds to a symbolic dynamics on two symbols.
\end{theorem}

Deformations of flows in general need not behave nicely with respect to the dynamics, and in particular a deformation can destroy the hyperbolicity and eliminate any number of periodic orbits for a flow. Our second result is that the deformation of the extended geometric model to the Lorenz equations cannot destroy any periodic orbits, but only continuously deform them, and thus we know that any knot types that appear as orbits for the geometric model appear also for the original equations. This is done using the cross section in order to apply methods from two dimensional dynamics \cite{Thurston1988diffeomorphisms, Boyland1994topologicalmethods, AsimovFranks1983Unremovable}.

\begin{theorem}\label{thm: periodic orbit isotopic}
Any knot in $S^3$ that is a periodic orbit in the extended geometric model $X_1$ is also a periodic orbit for the Lorenz equations at a trefoil parameter.
\end{theorem}

An immediate corollary of our results above and Ghys' theorem \cite{Ghys:KnotsDynamics} is thus the following.

\begin{corollary}\label{cor:modular}
Any knot that is a closed geodesic on the modular surface is also a periodic orbit of the Lorenz equations for some parameter.
\end{corollary}

Note that at the classical parameters it follows from Ghys' and Tucker's work that any knot appearing as a periodic orbit for the Lorenz equation is also a  closed geodesic on the modular surface, however, not every closed geodesic appears. Thus our theorem is in a sense a converse to Ghys'.

\subsection*{Acknowledgments}
The author wishes to thank Lilya and Misha Lyubich, Omri Sarig, Pierre Dehornoy, Sheldon Newhouse, Andrey Shilnikov, Eran Igra and the annonymous referees for helpful comments, And Amos Nevo and Joan Birman for their continuous encouragement.

\section{The existence of a heteroclinic trefoil}
In this section we prove the existence of a trefoil heteroclinic connection passing through the three singular points and infinity in the Lorenz equations for some parameter.
The first step in the proof is establishing the existence of a global cross section, and the second is an analysis based on the known homoclinic connections in the system.

The cross section we use is modified from the one commonly used in numerical studies which has a constant $z$ coordinate, and is reminiscent of Sparrow's study of the local  $z$ maxima along orbits \cite{Sparrow}. It remains a section also at parameter values for which the regular section fails to be one.

\begin{proposition}
There exists an open simply connected domain $A$ of parameters, where for any parameter  $(\beta,\sigma,\rho)$ there exists a two dimensional topological rectangle $R\subset \mathbb{R}^3$ with interior transverse to the Lorenz flow, so that the forward orbit of any point in $\mathbb{R}^3$ that does not limit onto one of the singular points meets $R$. 
\end{proposition}

\begin{proof}
We start by considering the hyperbolic paraboloid
\[
\mathcal{P}:=\{xy=\beta z\}
\]
Which is the set of points for which $\dot{z}=0$, containing the three singular points.
The paraboloid $\mathcal{P}$ divides $\R^3$ into two components, we call the part of $\R^3\setminus\PP$ above the origin the inside and the other component the outside.

For any orbit, its $z$ coordinate decreases while it is on the inside, and increases if and only if it crosses $\PP$ to be on its outside and so on.

A normal to $\PP$ is given by 
\[
N=(y,x,-\beta).\]
The vector field is tangent to $\PP$ exactly when it is orthogonal to the normal. Taking the inner product of $N$ with the vector field $X=(\dot{x}, \dot{y}, \dot{z})$ along the paraboloid $\PP$ we obtain
\begin{equation}\label{eqn:normal}
N\cdot X = \sigma y^2-\sigma xy+\rho x^2-xy-\frac{1}{\beta}x^3y
\end{equation}

The equation $N\cdot X=0$ is quadratic in $y$. There is always a solution at the origin and other solutions are given by the quadratic formula, where the discriminant is $(\sigma+1+\frac{x^2}{\beta})^2-4\sigma\rho$. Thus for any given $\sigma$ and $\rho$, there are always other solution curves once $|x|$ is large enough. If we want these solutions to be separated from the solution at the origin we need to require that $(\sigma+1)^2-4\sigma\rho<0$ so that the discriminant becomes positive only starting at some $|x_*|>0$.

Thus we will restrict our attention from now on to the domain 
\[
{A}=\left\{\beta,\sigma,\rho |\, \beta,\sigma>0, \rho>1, \rho>\frac{(\sigma+1)^2}{4\sigma}\right\}.
\]

For any parameter in $A$,
the set of points $N\cdot X=0$ for which the vector field is tangent to $\PP$ 
consists of a point at the origin, and of two one dimensional curves which we denote by $\delta_-$ and $\delta_+$ each containing one of the wing centers. the orbits cross $\PP$ to its outside below the two curves $\delta_-$ and $\delta_+$ and enter the inside above these curves. 
An example of the regions of entrance to the inside of $\PP$ and exit from it are depicted in Figure \ref{fig:paraboloid}.

\begin{figure}[ht]
    \centering
    \includegraphics[width=10cm]{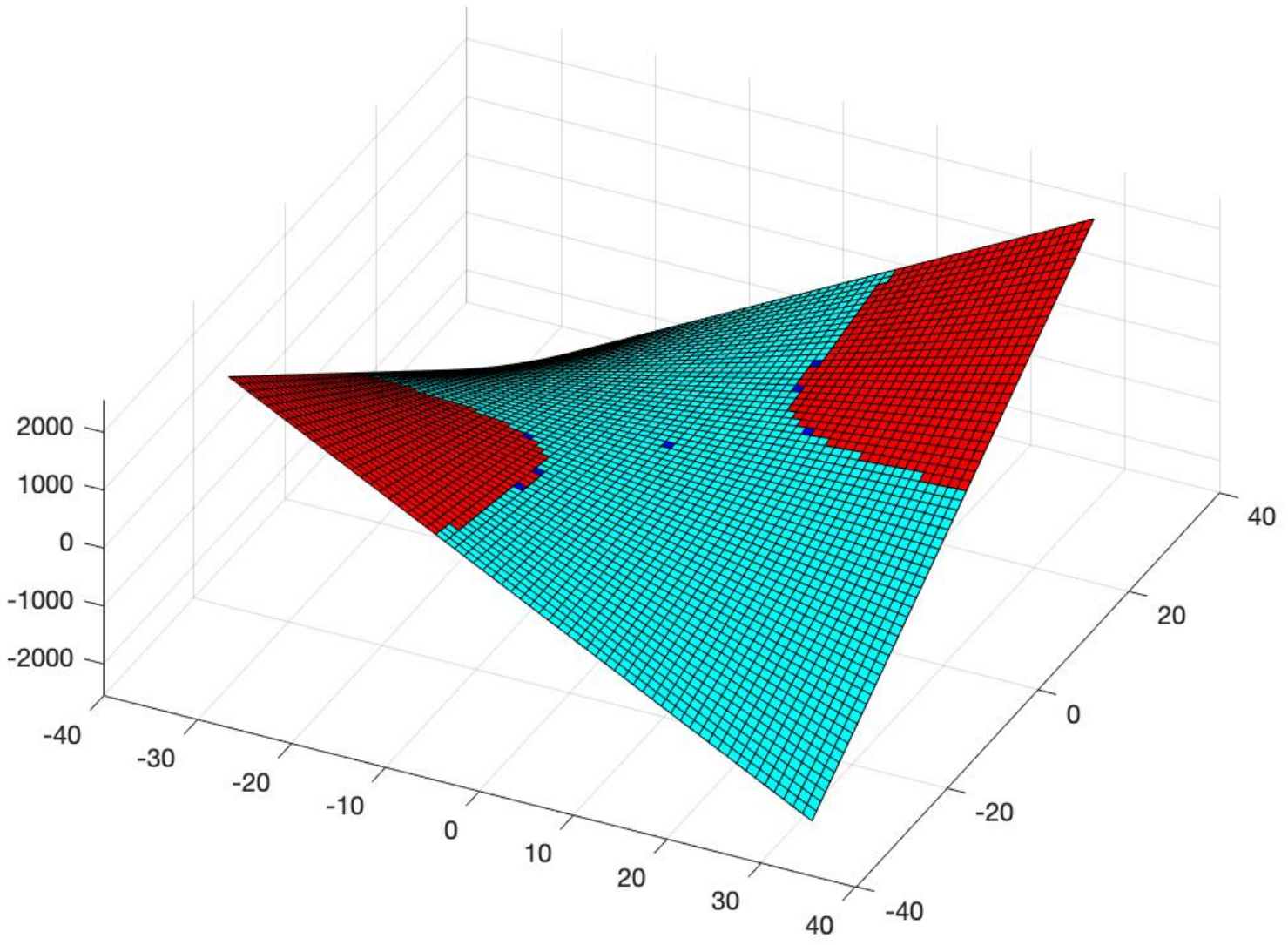}
    \caption{The paraboloid $\PP$ for the classical parameters, with the regions where the flow is directed inwards (toward the region above the origin) depicted in red (or dark) and outwards flow is depicted in blue.}
    \label{fig:paraboloid}
\end{figure}

Recall that an ellipsoid around the origin of the form
\[
\{\rho x^2+\sigma y^2+\sigma(z-2\rho)^2=C_0^2\}
\]
for a large enough $C_0$ is transverse to the vector field $X$, so that orbits only enter the region bounded by the ellipsoid (this is developed in  \cite{lorenz1963deterministic}. See page 197 of \cite{Sparrow} for estimates of $C_0$ for any parameter value).  We choose one of these ellipsoids transverse to the flow as a (topological) sphere around infinity and denote it by $S$. 

The Lorenz equation have a symmetry of rotation about the $z$ axis $\pi:(x,y,z)\mapsto(-x,-y,z)$, and $\pi(S)=S$.

The curves $\delta_-$ and $\delta_+$ will serve as two edges for the rectangle $R$. Note that they intersect the sphere $S$ each in two points, and $\pi(\delta_-)=\delta_+$.
We next define another pair of arcs, $\alpha_1$ and $\alpha_2$.
First, fix an $\varepsilon>0$. Denote the two intersection points of $\delta_\pm$ with $S$ by $(x_\pm^1, y_\pm^1, z_\pm^1)$  and $(x_\pm^2, y_\pm^2, z_\pm^2)$. Define four arcs $\omega_\pm^i\subset \PP\cap S$ continuing from each intersection point of $\delta_+\cap S$ downwards, until it reaches the plane $z=\varepsilon$. 

Define an arc $\omega^3$ connecting the endpoint of $\omega_+^1$ and $\omega_-^2$ in the plane $\{z=\varepsilon\}$, through one of the $xy<0$ quadrants, so that $\omega^3$ lies in the intersection of this plane and the ellipsoid, and in the inside of $\PP$. 

Denote the union $\omega_+^1\cup\omega_-^2\cup\omega^3$ by $\alpha_1$. Define $\alpha_2$ similarly by considering the images under the symmetry $\pi$ of $\omega_+^1,\omega_-^2,\omega^3$ and letting $\alpha_2$ be their union.        
The four arcs $\delta_-$, $\delta_+$, $\alpha_1$, $\alpha_2$ bound a closed topological rectangle $R$, contained in the union of $\PP$ and the plane $z=\varepsilon$. 

By construction, $R$ is preserved by the symmetry, $\pi(\alpha_1)=\alpha_2$ and $\pi(\delta_-)=\delta_+$.

Choosing $\varepsilon$ to be small enough, the part of $R$ that is contained in $\PP$ is disjoint from the set $N\cdot X=0$, and therefore the flow is transverse to $R$ on this set. The other part of $R$ is contained in the intersection of the plane $z=\varepsilon$ with the inside of $\PP$. Thus, the flow there satisfies $\dot{z}<0$, and as this part of the rectangle is horizontal, the flow is transverse to $R$ there as well. By smoothening $R$ in a small neighbourhood of the intersection of $\PP$ with the plane, we make interior of $R$ into a differentiable surface, and $R$ is transverse to the flow on all its interior.

Orbits cannot escape to infinity as the flow is pointing inwards through $S$. As $\PP$ contains no closed curve with constant $z$ coordinate, it follows that there are no nonsingular orbits with constant $z$ coordinate, the forward orbit of any point $x\in\R^3$ either has an increasing and then decreasing $z$ coordinate, in which case it must meet $R$, or else the orbit limits onto a singular point (or is fixed at a singular point), as required.
\end{proof}

\begin{lemma}
The component of the intersection of the stable manifold of the origin and $R$ containing the point on the $z$ axis  is a one dimensional curve $\eta$ dividing $R$ into two components, yielding a symbolic coding on two symbols.
\end{lemma}

\begin{proof}
Consider the linearization of the Lorenz flow,
\begin{equation}
    \begin{pmatrix}
    -\sigma & \sigma & 0 \\
    \rho-z & -1 & -x \\
    y & x & -\beta
    \end{pmatrix}.
\end{equation}
Using this, one finds that the linearization at the origin has one positive eigenvalue and two negative eigenvalues when $\beta>0$, $\sigma>0$ and $\rho>1$,  thus the stable manifold at the origin is two dimensional.  It is easy to see that the stable manifold contains the $z$ axis. Thus, the stable manifold intersects $R$ at least at one point, i.e. the point at which the $z$ axis intersects the plane $z=\varepsilon$.
From the transversality of the flow to $R$, and since the flow is tangent to the stable manifold, it follows that the intersection of the stable manifold and $R$ is transversal and thus the intersection is a one dimensional manifold $\eta$. The stable manifold is topologically a disk. It cannot have a periodic orbit on its boundary as such an orbit would link nontrivially the $z$-axis, in contradiction. Thus, the curve of intersection $\eta$ is a one-dimensional manifold that connects to infinity on both sides.

A small segment of the curve about the $z$ axis is contained in the plane $z=\varepsilon$, and consists of points whose orbits limit onto the origin without intersecting $R$ a second time. Orbits of points on $R\cap \PP$ flow outwards through $\PP$ then upwards, and to the inside of $\PP$ above $\delta_+$ or $\delta_-$, hitting  $R$ at least one more time before reaching a small neighbourhood of the origin. Thus, $\PP\cap\,\eta=\emptyset$, and we conclude that $\eta\cap R$ is contained in the plane $z=\varepsilon$. 

A priori it is possible that a disk within $\mathcal{W}^s(0)$ ripples through the intersection of the plane with the sphere $S$. In this case, choosing two points $q_1,q_2$ outside the sphere and an arc within $\eta$ connecting them, each $q_i$ is part of an orbit coming from infinity, and thus (as  $\mathcal{W}^s(0)$ is simply connected), all points  on the connecting arc are also part of orbits originating from infinity. Hence, we may modify $R$ to exclude the half disk bounded by the sphere $S$ and the arc connecting $q_1$ to $q_2$. Since orbits within the sphere cannot intersects the orbits connecting infinity to the subarc of $\eta$, all points in this half disk are on orbits arriving from infinity, and modifying $R$ will not alter its property that all orbits will hit it in forward time unless they limit onto the singular points.

It follows that the curve $\eta$ divides the (modified) rectangle $R$ into two symmetric subsets, each containing one of the wing centers. 
Denote the Poincar\'e return map by $\varphi$.
The curve $\eta$ yields a natural way to endow the system $(R,\varphi)$ with a symbolic dynamics with  two symbols, by recording in which of the two components of $R\setminus\eta$ is each point, and its subsequent images under $\varphi$.
\end{proof}

Denote by $A_1$ and $A_2$ the two components of $R\setminus\eta$.
Each component has one of the wing centers on its boundary, and thus the transition $A_1$ to $A_1$ appears in the dynamics, and also $A_2$ to $A_2$. Furthermore, below the plane $z=\varepsilon$ the manifold $\mathcal{W}^s(0)$ is fairly simple, in contrast to what happens above this plane (as is indeed observed numerically, e.g. in \cite{OSINGA2002VisualizingLorenz}). Points off $\eta$ return to $R$ without intersecting $\mathcal{W}^s(0)$, and thus orbits emanating from $A_1$ must enter $\PP$ through the quadrant $\{x,y>0\}$, and orbits from $A_2$ enter $\PP$ from the negative quadrant. Furthermore, points close to the entire edge of $A_i$ along $\eta$ are mapped close to a single point that is the first point where the seperatrix hits $R$. In this sense the return map send each of $A_i$ to a triangle like region and the dinamics is influenced by the location of the ``triangle tip". However, there is still too much freedom for this symbolic dynamics to be useful. 
This is why we next focus on special parameter values allowing us to pinpoint the behaviour of the return map.


\subsection{Existence of the heteroclinic trefoil}
\begin{proposition}\label{pro: half Stable to Infinity}
For any parameter in the domain $A$, one side of its stable manifold of each of the wing centers $p^\pm$ connects it to infinity.
\end{proposition}

\begin{proof}
Consider the half space
\[
\mathcal{Q}_+=\{(x,y,z)\,|\, x\geq x_0\},
\]
where $x_0=y_0=\sqrt{\beta(\rho-1)}$ and $z=\rho-1$ are the coordinates of the fixed point $p^+$. By considering the linearization (3) at $p^+$, plugging in its coordinates, 
one finds that there is always at least one negative real eigenvalue $\lambda$
and a corresponding eigenvector $v$ will have a first coordinate $v_1\neq0$. 

Next, consider $\partial\mathcal{Q}_+$. The flow is transverse to the half plane $\{x\equiv x_0, \, y<x\}$ which we denote by $\partial\mathcal{Q}_{out}$, pointing outwards at each point, and transverse and pointing into $\mathcal{Q}_+$ on  $\partial\mathcal{Q}_{in}=\{x\equiv x_0,\, y>x\}$.

We next consider the component of intersection of $R$ with $\mathcal{Q}_+$ containing $p^+$ which is a ledge $L$ sticking into $\mathcal{Q}_+$, with its interior transverse to the flow: It has one boundary that is a curve on the half plane $\{x\equiv x_0, \, y<x\}$ going down along its intersection with $\PP$ and continuing horizontally through the plane $\{z=\varepsilon\}$, and a second boundary curve on the half of $\delta_+$ within $\mathcal{Q}_+$ (note that these two boundary curves meet at $p^+$). By the definition of $R$, the ledge is a union of the planar region $\{x\geq x_0, y\leq 0, z=0\}$, and a region contained in the paraboloid $\PP$, smoothened along their intersection. The flow is transverse to the ledge, pointing downwards through the plane and outwards through $\PP$.

Consider next the half of $\delta_+$ inside $\mathcal{Q}_+$ and the orbits of all points on it under the flow.
Taking into account the equation defining $\delta_+$, it is monotone decreasing in $z$ when emanating from $p^+$, has a global minimum and is then monotone increasing. The flowlines starting on different points of $\delta_+$ enter immediately the inside of $\PP$ on the portion it is increasing. On the initial segment of $\delta_+$ where it is decreasing the flowlines exit $\PP$ at first. However, then their $z$ coordinate increases and they cannot exit $\mathcal{Q}_+$ through the half plane $\{x\equiv x_0, \, y<x\}$ without passing through $\PP$ first. They also cannot hit the other component of $R$ in $\mathcal{Q}_+$ by the continuity of the flowlines starting on this half of $\delta^+$, as all flowlines hitting the boundary of the second component along $\delta^+$ for $y>x$ come from inside $\PP$ (as it is monotone increasing in $z$). Thus, orbits emanating from the initial $\delta_+$ segment enter $\PP$ joining smoothly the orbits from the rest of $\delta_+$.
All orbits starting at $\delta_+$ must at this point either wind around $\delta_+$ and intersect the ledge $L$, or exit $\mathcal{Q}_+$ through the part of $\partial\mathcal{Q}_{out}$ inside $\PP$.
In both cases, the union of the shell of flowlines $\bigcup_{t\geq0}\psi^t(\delta_+)$ together with the ledge and the outflowing boundary $\partial\mathcal{Q}_{out}$ bounds  an open cone with tip at $p^+$. At any point on the boundary of this cone,  the flow is either tangent to it or points outwards from it. Therefore,  a component of the stable manifold of $p^+$ must be contained in this cone. 

Considering the transverse sphere $S$, the cone is foliated by its images along the flow, which are parallel disks. Thus the arc of the stable manifold of $p^+$ in the cone connects $p^+$ to $\infty$ trivially in $\mathcal{Q}_+$.
The same claim for $p^-$ follows from the symmetry.
\end{proof}

This has been observed numerically as well, see for example \cite{Sparrow} and \cite{creaser15alpha}.

\begin{theorem}\label{thm: Trefoil existence}
There exists a point in the parameter space at which there exists a heteroclinic trefoil for the Lorenz equations.
\end{theorem}

A central tool of this proof is the following theorem, proving existence of homoclinic orbits and their linking with the vertical lines
\[
L^\pm=\{x=y=\pm\sqrt{\beta(\rho-1)}\}
\]
through the two wing centers.

\begin{theorem}[Chen, Theorem 1.1 and Lemma 4.3 of \cite{Chen1996RandomlyRotatedHomoclinic}] \label{thm:chen}
For any given positive number $\beta$ and non-negative integer $N$, there exists a large positive constant $\rho_0(N,\beta)$ such that for each $\rho>\rho_0$ there is a positive number $\sigma=\sigma(N,\beta,\rho)=(2\beta+1)/3+O(\rho^{-1/2})$ such that the Lorenz system has a homoclinic orbit associated with the origin which rotates around $L^+$ exactly  $\frac{N+1}{2}$ times, and rotates $\frac{N-1}{2}$ times around $L^-$.
\end{theorem}

\begin{proof}[proof of Theorem~\ref{thm: Trefoil existence}]
The description in Theorem~\ref{thm:chen} of the homoclinic orbits implies that the positive half of the  simplest such orbit corresponding to $N=1$ links once with $L^+$ and does not link with $L^-$. It follows that the separatrix, i.e. the unstable manifold of the origin, at that parameter returns to $\PP$ from the $x,y>0$ quadrant, linking $\delta_+$, and then connects to the origin after hitting $R$ a single time.
At the second homoclinic loop $N=2$ the separatrix links $\delta_+$ like in the first loop, and then continues to also link once with $\delta_-$. It then connects to the origin after it hits $R$ the second time.

The main tool in Chen's work is a change of variables that changes the Lorenz equations to the system
\begin{equation}\label{sys:chen}\left\{\begin{array}{rlr}
\dot{u}  & = v\\
\dot{v} & = u-2u^3-2\varepsilon(uw+v+a\varepsilon w) \\
\dot{w} & = v^2-\varepsilon b w
\end{array}\right.
\end{equation}
Where $a, b$ and $\varepsilon$ are real parameters depending on $\beta, \sigma,\rho$, and $\varepsilon$ approaching 0 when $\rho$ approaches $\infty$. This change of variables shows that when $\rho$ is large the system is a small perturbation of a system that is a Duffing equation for the first two variables with two sinks and one saddle point at the origin. It thus follows that increasing $\rho$ along the curve in parameter space with fixed $\beta$ corresponding to the first homoclinic connection in Chen's theorem that the homoclinic orbit is on the boundary of the basin of attraction of the wing center it encircles.

Consider a path $\kappa(s)$ connecting two points in parameter space $\kappa(0)$ corresponding to the first homoclinic loop $N=1$ and $\kappa(1)$ to the second $N=2$.
By Chen's result we can choose the entire path $\kappa(s)$ to lie in a region $\rho>M$ for any fixed $M$. We can thus assume that the orbits of the system \ref{sys:chen} with $\varepsilon=0$ approximates the orbits of the Lorenz system along $\kappa(s)$ as closely as we wish.
We may assume without loss of generality that $\kappa(s)$ is not a parameter on the first homoclinic curve for $s>0$ (that is $\kappa(s)$ is transverse to the homoclinic curve and does not return to it). From the linking of the second homoclinic with $L^\pm$ we know that at $k(1)$ the tip of the right half $A_1$ of the cross section $R$ is mapped to the left (and then hits the stable manifold along $\eta$ at its second return), it follows that the return map $r(s)$ for a small positive $s$ along $\kappa$ takes the triangle tip corresponding to the right half of the cross section to the left of $\eta$ (into $A_2$), and the tip of the left part $A_2$ of $R$ to the left.

As explained heuristically in \cite{Sparrow}, and rigorously worked out in \cite{chen1995HyperbolicSets}, for a small positive $s$ as the triangle tips are close to $\eta$ and the flow is expanding in a direction transverse to $\eta$ in a small enough neighborhood, the subsequent images of the triangle tips get further away from $\eta$ (this is the basis for the ``homoclinic explosion", i.e. the creation of infinitely many periodic orbits in a thin strip near $\eta$ for a small $s>0$). 
Thus, the orbits of the seperatrix, including the triangle tip enter the basin of attraction of the wing centers (as this basin is stable under perturbation) and therefore limits onto the two wing centers for some parameter range $k(s)$ for $s$ positive but small enough.

Let $\kappa_*$ be such a heteroclinic parameter. 
It follows from the proof that the positive half of the separatrix winds once around $\delta_+$ and connects to $p^-$, and the other half of the separatrix connects to $p^+$ on the other side of the origin.
It now follows from \ref{pro: half Stable to Infinity} that the continuation of the connection to the other side of each wing center along its stable manifold connects trivially to infinity in $\mathcal{Q}_+$ and its image under $\pi$, and therefore the resulting knot is a trefoil knot as required.
\end{proof}

\begin{remark} 
The existence of a first order heteroclinic connection in an open domain near the first homoclinic curve is known numerically, see for example Figure 3 in \cite{barrio12knead}, originally appearing in \cite{shilnikov1980bifurcation}. In fact it seems numerically that the trefoil heteroclinic connection exists in an open domain to the right of the entire homoclinic curve and not just for large $\rho$. This heteroclinic connection is not a connection that experts call a T-points which is always unstable, existing on a codimension 2 set of parameters. An example for a parameter value at which the stable trefoil heteroclinic connection exists is $\beta=8/3$, $\sigma=10$ and $\rho=17$, shown in Figure~\ref{fig:stable} below.
\end{remark}

\begin{figure}[ht!]
\centering
\begin{overpic}[width=9cm]{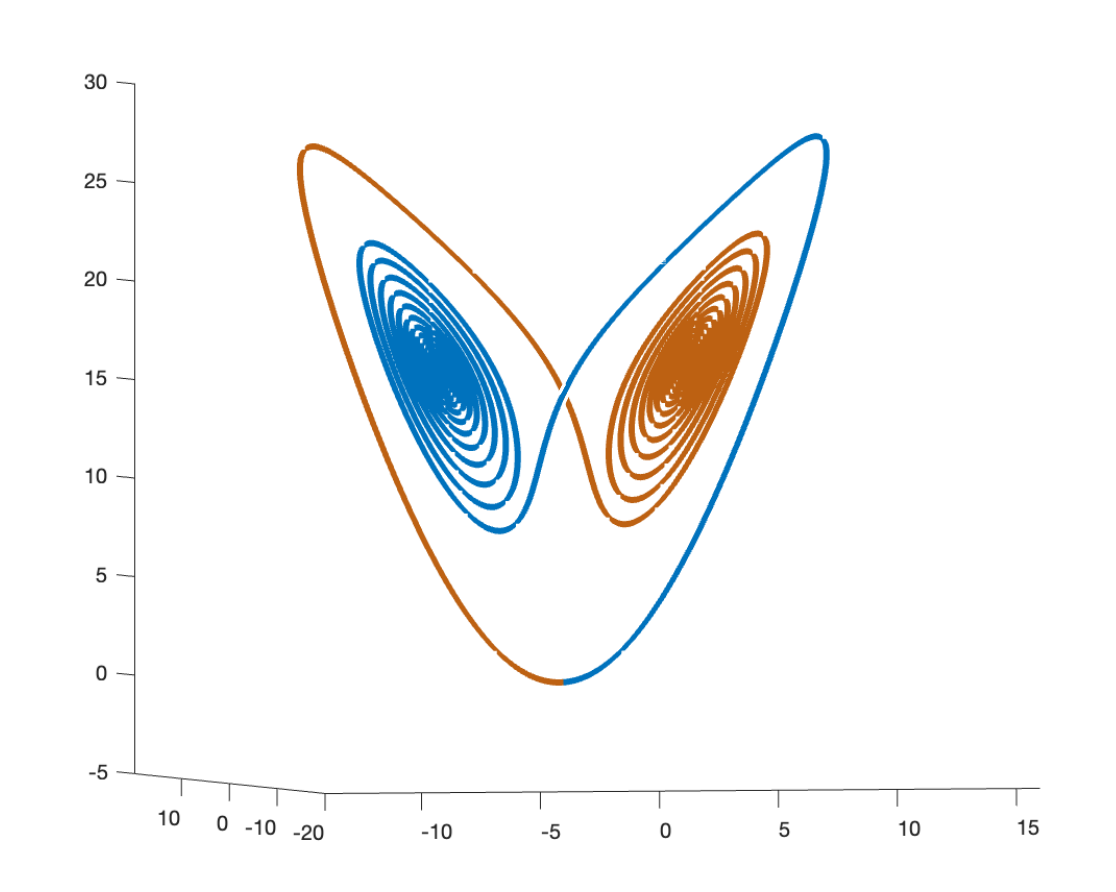}
    \put(1,72){$z$}
    \put(90,0){$x$}
    \end{overpic}
\caption{A pair of stable heteroclinic connections that can be completed to a trefoil knot, depicted for $\beta=8/3$, $\sigma=10$ and $\rho=17$.}\label{fig:stable}
\end{figure}

\begin{remark}
Theorem~\ref{thm: Trefoil existence} proves that a parameter for which there is a heteroclinic trefoil exists at any fixed $\beta$, as Chen's theorem applies and the path $\kappa(s)$ can be chosen to lie in any plane with fixed $\beta$ value.
\end{remark}

\section{The dynamics at a trefoil parameter}

In this section we prove Theorems \ref{thm:MainThm} and \ref{thm: periodic orbit isotopic}. 

\begin{theorem}\label{thm:isotopy}
At any parameter value where the Lorenz flow has a cross section and a trefoil heteroclinic connection, the flow defined by the Lorenz equations is connected through a one parameter family of flows, to the extended geometric model $X_1$ of \cite{BonattiPinsky2021}. For any flow in the one parametric family there exists a global
cross section, so that the return map corresponds to a symbolic dynamics on two symbols.
\end{theorem}

\begin{proof}
Consider a parameter $\kappa$ where the  heteroclinic connection is a trefoil knot, and assume that the cross section $R$ exists at this parameter. As before, $R$ is divided into two components $A_1$ and $A_2$ by the stable manifold of the origin, and let $\phi$ be the first return map to $A_1\cup A_2$ under the flow.

The image of $A_2$ (and likewise of $A_1$) under $\phi$ is a topological disk that includes both wing centers on its boundary: It includes the wing center $p^-$ that is on the boundary of $A_2$ as the fixed point is is equal to its image, and it also must include the wing center $p^+$ because a sequence of points in $A_2$ limiting onto $\eta$  have flowlines that continue near the stable manifold of the origin, enter smaller and smaller neighborhoods of the origin and then continue along the separatrix (which is part of the trefoil) that connects to $p^+$ (see Figure~\ref{fig:connections}). Thus for any trefoil parameter the image of $A_i$ covers both $A_1$ and $A_2$ for the original equations.

\begin{figure}[ht!]
\centering
\begin{overpic}[width=4.5cm]{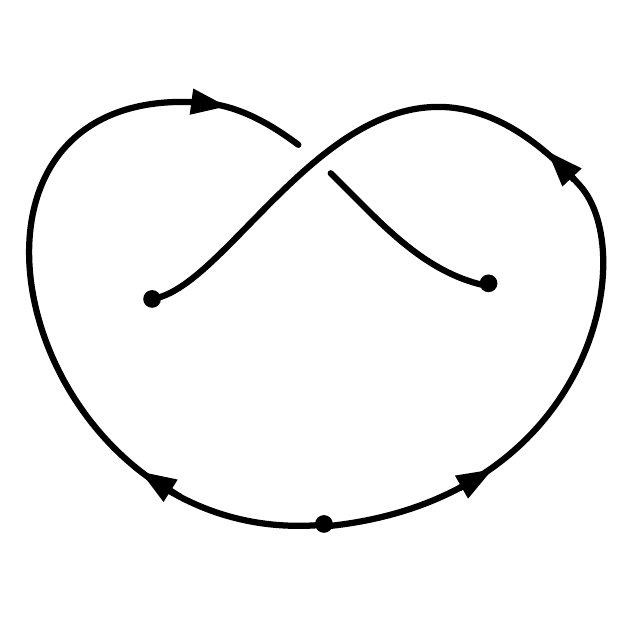}
\put(30,48){$p^-$}
\put(80,48){$p^+$}
\put(50,2){$0$}
\end{overpic}
\caption{The two seperatrices at a trefoil connection paramater $\kappa$.}\label{fig:connections}
\end{figure}

We now turn to define the deformation of the flow. As we prove the theorem for any parameter where we have a heteroclinic trefoil, the connection might not be stable as given by Theorem \ref{thm: Trefoil existence}, but could be unstable as in the parameter of the first T-point (where the wing centers are saddle-foci and not sinks). 
In this case, the first part of the deformation turns the trefoil into a stable one. Then we use a trick employed also in \cite{BonattiPinsky2021} (c.f. \cite{Ghys:KnotsDynamics}), which is blowing up the one dimensional trefoil into a thick three dimensional trefoil, which allows us to continue and deform the flow further in the compactified trefoil complement.

Step 1: In the case the two wing centers are saddles having a two dimensional unstable manifold we deform the flow, introducing a \emph{subcritical Hopf bifurcation} turning the wing centers into sinks. This is done by adding along the two dimensional unstable manifold a component to the vector field tangent to the manifold and pointing toward the wing center. This adds a trivial periodic orbit encircling each wing center, corresponding to an additional saddle fixed point in each of the rectangles $A_1$ and $A_2$. The deformation has a natural parameter, so that the vector field $X_\kappa$ for $\kappa$ a trefoil parameter in $A$ is replaced by a deformed field $X_{\kappa,r}$ where $r$ is the radius in which the deformation is preformed, and can be as small as required.

Note that for the deformed field $X_{\kappa.r}$ points that are very close to $\eta$ flow very close to the stable manifold reaching a neighbourhood of the origin, and continue very close to the seperatrices reaching the basin of attraction of either $p^+$ or $p^-$. Thus there is a neighborhood of $\eta$ that limits onto either one of the wing centers.

\begin{figure}[ht!]
\centering
\begin{overpic}[width=4.5cm]{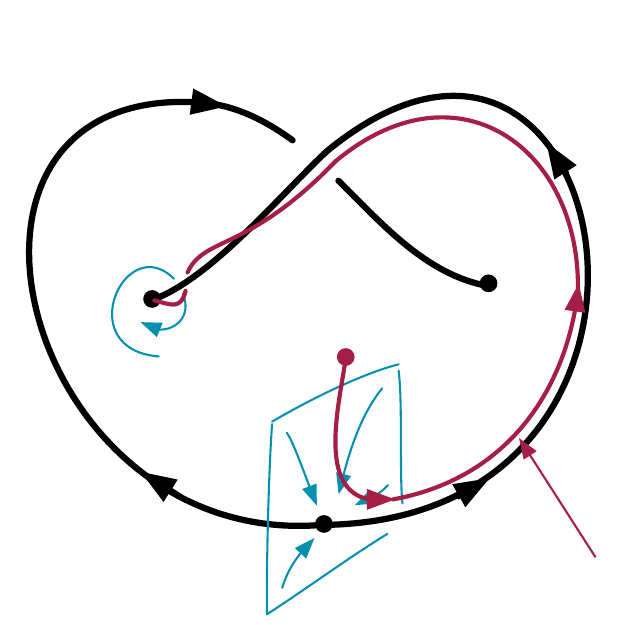}
\put(30,48){$p^-$}
\put(80,48){$p^+$}
\put(73,7){orbit starting near $\eta$}
\end{overpic}
\caption{After the deformation, an orbit close enough to $\eta$ will limit onto one of the wing centers (where ``close enough" will depend on the deformation parameter $r$).}\label{fig:connections2}
\end{figure} 

In the case the two wing centers are sinks to begin with, Step 1 is redundant as we are already in the case where there exists at least one additional fixed point in each of $A_1$ and $A_2$, and points close enough to $\eta$ limit onto one of the wing centers. 

\begin{figure}[ht!]
\centering
\begin{overpic}[width=6cm]{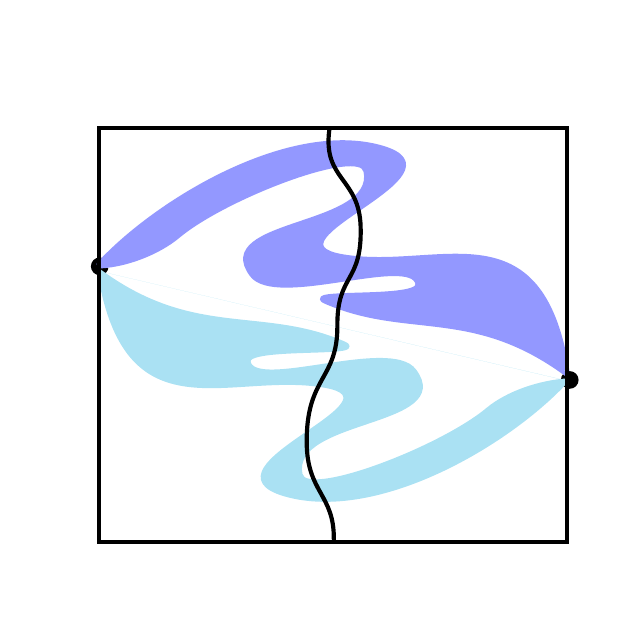}
\put(25,16){$A_1$}
\put(65,62){$\phi_\kappa(A_1)$}
\put(21,38){$\phi_\kappa(A_2)$}
\put(73,16){$A_2$}
\put(50,5){$\eta$}
\end{overpic}
\begin{overpic}[width=6cm]{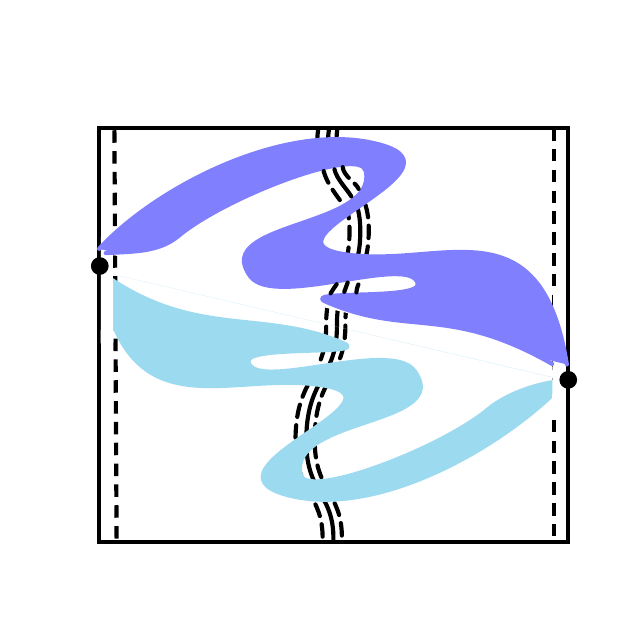}
\put(25,16){$A_1$}
\put(62,62){$\phi_{\kappa.r}(A_1)$}
\put(20,38){$\phi_{\kappa,r}(A_2)$}
\put(73,16){$A_2$}
\put(50,5){$\eta$}
\end{overpic}
\caption{An illustration of the action of the deformation on the return map to the cross-section.}\label{fig:return}
\end{figure}

Step 2: Remove a three dimensional neighborhood of the trefoil knot as in \cite{BonattiPinsky2021} in the following way: Start from a very small neighborhood of the origin and extend it along the separatrices, so that the flowlines point inwards all along it, until reaching the basin of attraction of the wing center on either side. There, let the neighborhood grow wider, so that its boundary connects to the trivial orbit encircling each wing center. Continue to widen the neighbourhood while continuing along the trefoil on the other side of the wing center, until it connects to the sphere around infinity (that it connects to by Proposition~\ref{pro: half Stable to Infinity}). Using the local form of the flow in the neighborhoods of the fixed points and their invariant manifolds, it is shown in \cite{BonattiPinsky2021} that the neighborhood can be chosen so that its boundary is a two dimensional torus, transverse to the flow generated by $X_{\kappa,r}$ except along two tangent orbits.
The action of this deformation on the cross-section is depicted in Figure~\ref{fig:return}, where as a neighborhood of the wing centers is contained in the trefoil neighborhood, flowlines starting close enough to $\eta$ exit the manifold into this neighborhood and never return to the cross section.

The third step of the proof is now to define the one parameter family taking the flow $X_{\kappa,r}$ to the extended geometric model, in which the cross section is uniformly stretched in one direction and contracted in another. In order to make use of methods from two dimensional dynamics, we embed the return map $\phi_{\kappa,r}$ into a diffeomorphism of a closed surface.

To this end consider the action of the Arnold cat map 
$\begin{pmatrix}
1 & 1 \\ 1 & 2
\end{pmatrix}
$
on the two dimensional torus $\mathbb{R}^2/\mathbb{Z}^2$. It has a fixed point at the origin and we can  arrange by an isotopy of the diffeomorphism to have a second fixed point. The diffeomorphism can naturally be defined on the complement of the fixed points, yielding a diffeomorphism on the twice punctured torus. For any punctured surface, one can choose a one dimensional graph so that the surface can be pushed into a small neighborhood of the graph. Such a graph is called a \emph{spine}, and any map of the spine to itself has a natural way to be extended to its neighborhood, i.e. to a diffeomorphism of the entire surface (in our case, of the punctured torus) \cite{BestvinaHandel1995TrainTracks, Boyland1994topologicalmethods}.

Choose a graph that is a spine for the punctured torus consisting of a vertical edge $w$, a diagonal edge $v$, and a small circle around the fixed point at the origin. The action of the map on an a neighborhood of the graph is given in Figure~\ref{fig:Arnold}, and we embed the action of $\phi_{\kappa,r}$ into the a map of the torus as the action on the two components of $v\cap \psi_C^{-1}(v)$. 

\begin{figure}[ht!]
\centering
\begin{overpic}[width=5cm]{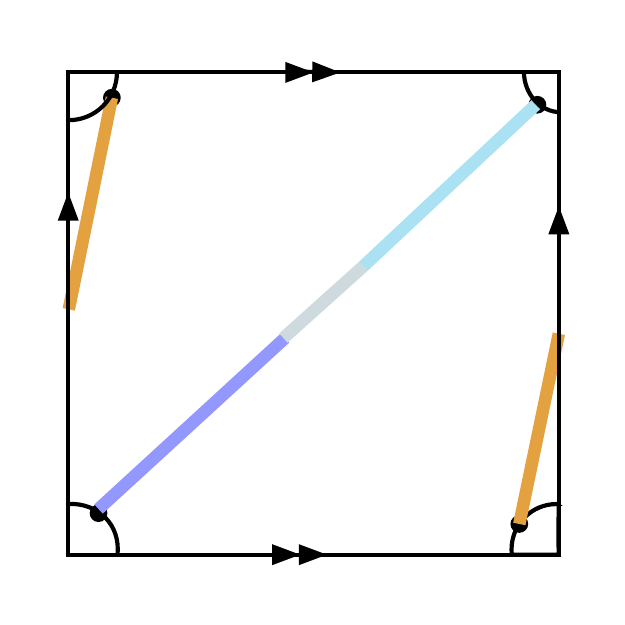}
\put(79,30){$w$}
\put(50,55){$v$}
\end{overpic}
\begin{overpic}[width=5cm]{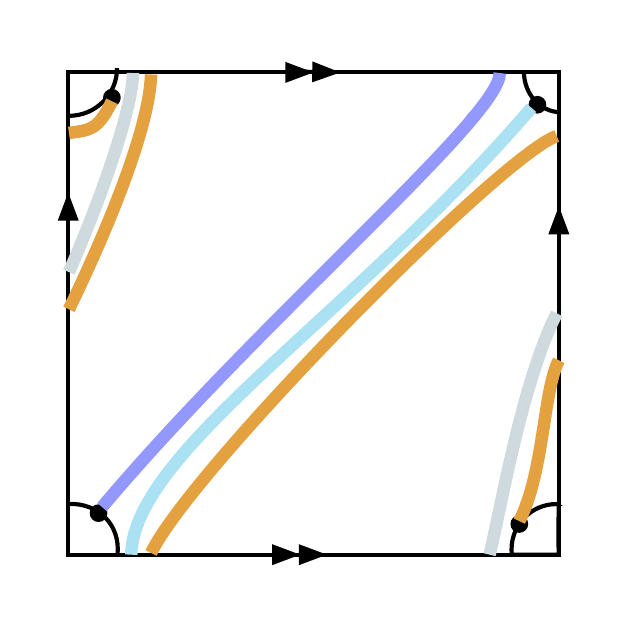}
\end{overpic}
\begin{overpic}[width=5cm]{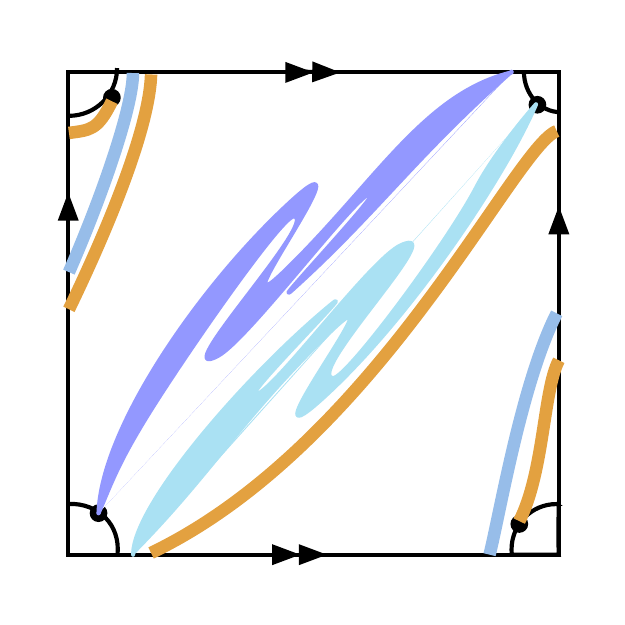}
\end{overpic}
\caption{The edges $v$ and $w$ and their images under the cat map, and the embedding of the action of $\phi_{\kappa,r}$ on $A_1$ and $A_2$ taken as neighborhoods of subarcs of $v$.}\label{fig:Arnold}
\end{figure}

Thus this procedure uses $\phi_{\kappa,r}$ to define a diffeomorphism $\Psi_{\kappa,r}$ of the once punctured torus to itself for any $r>0$. Any diffeomorphism of a surface is isotopic to a Thurston Nielsen canonical form \cite{Thurston1988diffeomorphisms} that is hyperbolic (in this context called pseudo-Anosov), finite order or reducible, and it is the most efficient in its isotopy class. There is an algorithmic way to obtain the canonical form called the Bestvina Handel algorithm \cite{BestvinaHandel1995TrainTracks}. Applying the algorithm in general requires several steps that modify the spine together with its graph map to get a map with smaller topological entropy, until it terminates at a graph map defining an isotopic  diffeomorphism that is the canonical form. In the case of the diffeomorphism $\Psi_{\kappa,r}$, applying a single step of the algorithm which is tightening the map along the diagonal edge results in the cat map, which (since it is hyperbolic) is already in the canonical form. Pulling the isotopy restricted to a neighborhood of $w$ back to the cross section (indeed the isotopy given by the Bestvina Handel algorithm is restricted to a neighborhood of the graph) defines an isotopy between $\psi_{\kappa,r}$ and a hyperbolic map of $R$ that stretches $A_1$ and $A_2$ uniformly, preserving both a horizontal unstable and a vertical stable foliation (given both by pullback of these foliations from the torus). This isotopy of the map from the cross section $R$ to itself induces an deformation of the flow by deforming continuously the flowlines to their new image under the map just before they intersect $R$.

\begin{figure}[ht]
    \centering
    \begin{overpic}[width=7cm]{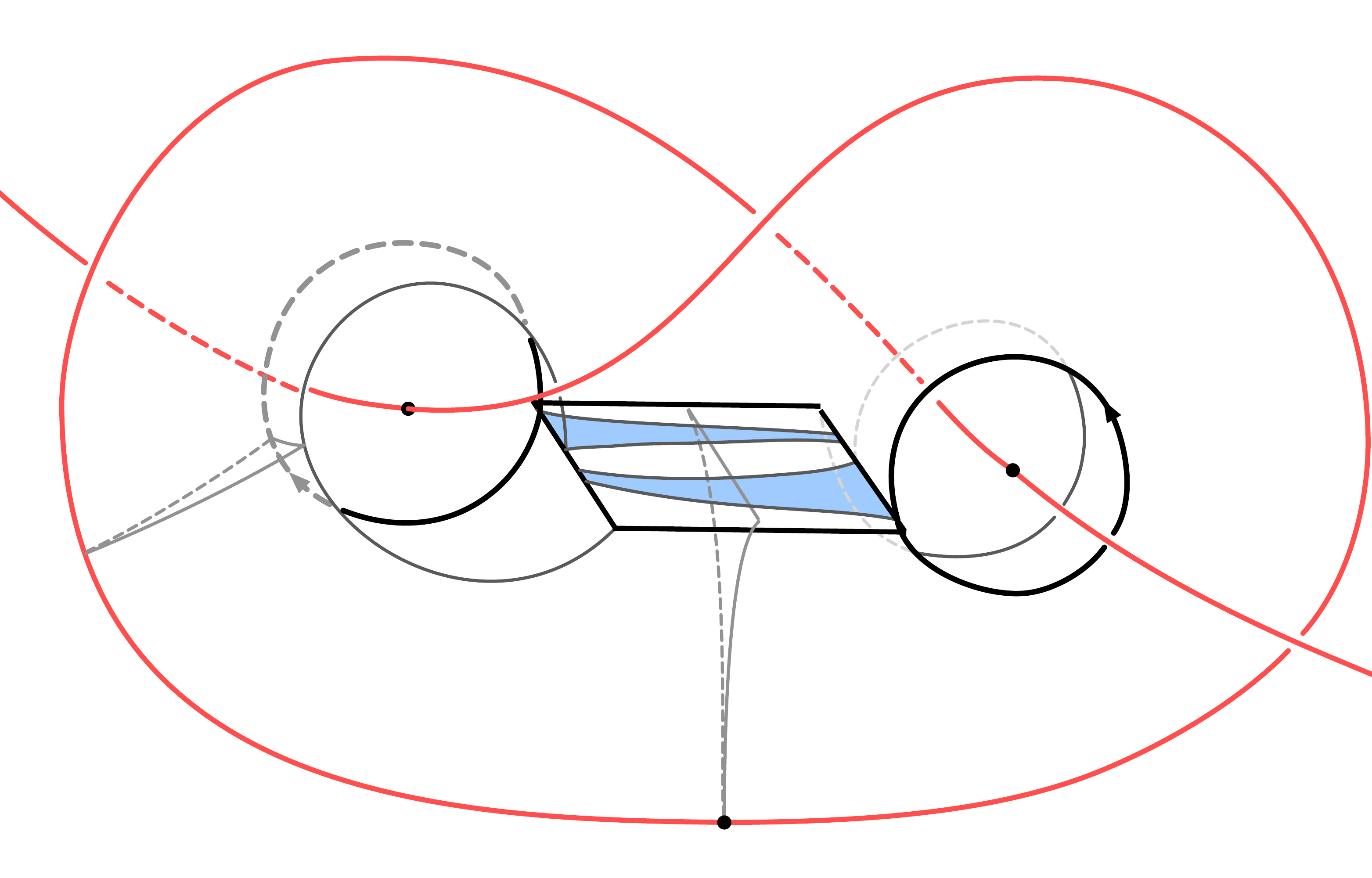}
    \put(30,30){$p^-$}
    \put(72,32){$p^+$}
    \put(54,-2){$0$}
    \end{overpic}
    \caption{The extended geometric model, defined on $S^3$ and containing an invariant trefoil \cite{BonattiPinsky2021}. A neighborhood of the trefoil can be removed and the resulting flow is the unique hyperbolic flow on the trefoil complement with the given two to one return map.}
    \label{fig:geometric}
\end{figure}

To conclude the proof, it follows from the theory of B\'eguin-Bonatti's geometric models (Theorem 0.3 of \cite{beguin2002flots}) that there is a unique flow on the trefoil complement with a hyperbolic two-to-one return map, which is the extended geometric model $X_1$ constructed in \cite{BonattiPinsky2021} and depicted in Figure~\ref{fig:geometric}. Thus, the one parameter family defined by concatenating the one in Step 1, increasing $r$ continuously from 0 to some $0<r_0$, with the one in Step 3 for this fixed $r=r_0$ yields a family of flows $X^s_\kappa$ connecting the Lorenz flow $X_{\kappa}=X_\kappa^0$  to the extended geometric model $X_\kappa^1$, as required. The exact value of $r_0$ is not important, as the resulting hyperbolic flows are all topologically  equivalent due to the structural stability.
\end{proof}

The proof of Theorem~\ref{thm:MainThm} now follows from Theorems~\ref{thm: Trefoil existence} and \ref{thm:isotopy}.

\begin{proof}[proof of Theorem~\ref{thm: periodic orbit isotopic}]
Consider the one parametric family $X^s_\kappa$ given by Theorem~\ref{thm:isotopy}, deforming the Lorenz flow at a trefoil parameter $\kappa$ to the hyperbolic extended geometric model. At $s=1$ the return map to $R$ embedded into the neighborhood of the diagonal of the torus is a restriction of the hyperbolic cat map $h_1$. For any hyperbolic diffeomorphism $f_0$ of a surface, each of its periodic orbits is \emph{isotopy stable} (see \cite{AsimovFranks1983Unremovable, BirmanKidwell1982FixedPoints} and Theorem 7.4 in \cite{Boyland1994topologicalmethods}). This means that for any orbit $x_0$ and any isotopy $f_t$ from $f_0$ to another diffeomorphism, a  periodic orbit $x_t$ exists for any $f_t$. Hence any periodic orbit for the cat map $h_1$ there exists an isotopic orbit for the return map $\phi_{\kappa,r}$ of $X^s_\kappa$, $s>0$. Thus, as the isotopy $X_\kappa^s$ for $s>0$ is a pullback of an isotopy on the two dimensional torus, all the periodic orbits of the extended geometric model continue to exist along the isotopy. This remains true when the deformation starts with a stable trefoil connection, as for any stable trefoil the return map can be embedded as a diffeomorphism of a punctured torus and thus the deformation corresponds to an isotopy of the punctured torus which preserves all these orbits.

Therefore, at this point, if we started with a parameter in the Lorenz equation for which the trefoil is stable we are done.
The last step in the proof is to show that the isotopic orbits exist also for $s=0$, i.e. when reducing the Hopf bifurcation parameter to zero. 
To this end, we use results from \cite{alligoodMalletYorke2006index} and \cite{ChowMallet-ParretYorke1983Index}, where an index can be defined for a periodic orbit, and if the index is nonzero the orbit persists under isotopy.
From the hyperbolicity on the trefoil complement of the extended geometric model it follows that its orbits all have non zero index \cite{ChowMallet-ParretYorke1983Index}, and have minimal period, and these properties persist along the isotopy for the isotopic orbits of the deformed flow for any $s>0$.
At $s=0$ the trivial Hopf orbits collapse onto the wing centers and disappear. All other points continue also for $s=0$, again by \cite{alligoodMalletYorke2006index}, as their indices are nonzero, they have minimal period, and they remain bounded away from the singular points. Thus, all periodic orbits except  possibly the Hopf orbits continue in the same isotopy class in the limit $s=0$, which is the original Lorenz system at the trefoil parameter we started with.
\end{proof}

\begin{remark}
Another family of flows having a return map given by the Arnold cat map was considered by Birman and Williams in \cite{BirmanWilliams2}. There, the family of flows are transverse not only to one cross section but to an entire two dimensional fibration of the space. Similar to our situation, the set of periodic orbits of any flow in this family contains, up to isotopy,  the set of all orbits corresponding to the cat map (which is remarkably all possible knots in $S^3$ \cite{ghrist1997AllLinks}).
\end{remark}

\begin{proof}[Proof of Corolary \ref{cor:modular}]
Every closed orbit on the Lorenz template, and therefore in the extended geometric model $X_1$ is, as a knot, a closed geodesic on the modular surface, except the two trivial orbits that are the template boundaries \cite{Ghys:KnotsDynamics}, or the trivial Hopf orbits in the extended model \cite{BonattiPinsky2021}.
By the proof of Theorem~\ref{thm: periodic orbit isotopic}, every such knot is also a periodic orbit for the Lorenz equation at a trefoil parameter, that exists by Theorem \ref{thm: Trefoil existence} (whether or not this parameter is beyond the Hopf bifurcation curve).
\end{proof}

\begin{remark}
Giving a full analytical solution to Smale's 14th problem will require a way to prove the nonexistence of attracting periodic orbits for the Lorenz equations, showing the $A_1,\,A_2$ symbolic partition is generating. This seems to be hard as it is a local phenomenon that cannot be obstructed topologically. Furthermore, it is observed numerically  that there are parameters for which such orbits do appear, but at $\rho$ values of about $100$ and above, while trefoil parameters appear for much lower $\rho$ values as well.
\end{remark}

\printbibliography
\end{document}